\documentclass[12pt]{amsart}

\usepackage{mathrsfs}
\input epsf
\usepackage{amssymb}
\usepackage{amsfonts}
\usepackage{amsmath}
\usepackage{mathrsfs}
\usepackage{amsthm}
\newtheorem{theorem}{Theorem}
\newtheorem{lemma}[theorem]{Lemma}

\theoremstyle{definition}

\theoremstyle{remark}

\newcommand{\C}{\mathbb{C}}
\newcommand{\R}{\mathbb{R}}

\begin{document}

\title[ Mahler's conjecture and box splines]{A remark about  Mahler's conjecture and the maximum value of box
splines }
\author{
 Zhiqiang Xu}

\begin{abstract}
In this paper, we recast a special case of Mahler'c conjecture by
the maximum value of box splines. This is the case of polytopes
with at most $2n+2$ facets. An asymptotic formula for univariate
box splines is given. Based on the formula, Mahler's conjecture is
proved in this case provided $n$ is big enough.

\end{abstract}

\maketitle

\section{introduction}

Let $K$ be a symmetric convex body in $\R^n$, and let $K^*$ be its
polar $\{x:|\left<y,x\right>|\leq 1 \mbox{ for all } y\in K\}$. An
old conjecture of  Mahler is
\begin{equation}\label{eq:mahler}
{\rm vol}(K)\cdot {\rm vol}(K^*)\,\,\geq\,\, \frac{4^n}{n!}.
\end{equation}
We note that an $n$-dimensional parallelepiped, which has $2n$
facets, gives equality in (\ref{eq:mahler}). So, the first
non-trivial case of Mahler's conjecture is the symmetric convex
body $K$ with $2n+2$ facets. Such $K$ can be realized, up to
affine invariance, as a one-codimensional section of an
$(n+1)$-dimensional cube. This case has been raised as a separate
problem by Ball \cite{mahler}. Ball also shows an interesting
relation between the special case of Mahler's conjecture and
solutions of a scaling equation. For each $r\in \R$, we set
\begin{equation}\label{eq:def}
\varphi_A(r)\,\,:=\,\,{\rm vol}\left((H+rA)\cap Q_n\right)
\end{equation}
where $Q_n:=[-\frac{1}{2},\frac{1}{2}]^n$, $A:=(a_1,\ldots,a_n)$
is an unit vector and $H:=\left<A\right>^{\bot}$. Throughout this
paper, without loss of generality, we suppose $0<a_1\leq a_2\leq
\cdots \leq a_n$.
  In particular, Ball shows that the special case of Mahler's
  conjecture is equivalent to
  \begin{equation}\label{eq:ba}
  \varphi_A(0)\cdot E\left(\left|\sum_{k=1}^n a_k
\varepsilon_k\right|\right)\,\,\geq\,\, 1,
  \end{equation}
where $\varepsilon_k$ is a sequence of mutually independent random
variables with distribution $P\{1\} \,=\,P\{-1\}\, =\,1/2$. In
\cite{lopez}, the authors prove (\ref{eq:ba}) for $n\leq 8$ by
reducing the problem to a search over a finite set of polytopes
for each fixed dimension.

An interesting observation is that $\varphi_A(\cdot)$ is a box
spline, a popular tool in approximation theory. Therefore, we can
recast the special case of Mahler's conjecture by the maximum
value of box splines. Using the saddle point approximation in
statistics, we give an asymptotic formula of $\varphi_A(\cdot)$.
Based on the asymptotic formula, we obtain the following result:
\begin{theorem}\label{th:main}
Suppose that there is a constant $c_0$ so that $a_n/a_1<c_0$ for
any $n$. Then there exists a positive integer $N_0(c_0)$ so that
$$
 \varphi_A(0)\,\cdot \, {E\left(\left|\sum_{k=1}^n
a_k \varepsilon_k\right|\right)}\,\,\, \geq\,\,\, 1
$$
when $n\geq N_0(c_0)$.
\end{theorem}

This paper is organized as follows. In Section 2, we introduce box
splines and show the relation between box splines and the special
case of Mahler's conjecture. In Section 3, we use the saddle point
approximation to give an asymptotic formula of univariate box
splines. Section 4 presents the proof of Theorem \ref{th:main}.

\section{Box splines}
Suppose that $M$ is a $s\times n$ matrix. The {\it box spline}
$B(\cdot |M)$ associated with $M$ is the distribution given by the
rule \cite{deboor1,deboor2}
\begin{equation}\label{Eq:boxspline}
\int_{{\R}^s}B(x |M)\phi (x)dx = \int _{[0,1)^n}\phi (Mu)du,\,\,\,
\phi \in {\mathscr D}({\R}^s).
\end{equation}
By taking $\phi={\exp}(-i\zeta\cdot )$ in (\ref{Eq:boxspline}), we
obtain the Fourier transform of $B(\cdot |M)$ as
\begin{equation}\label{eq:boxcube}
\widehat{B}(\zeta
|M)\,\,=\,\,\prod_{j=1}^n\frac{1-\exp(-i\zeta^Tm_j)}{i\zeta^Tm_j},\,\,\,\,\,
\zeta \in {\C}^s.
\end{equation}
The following formula shows the relation  between box splines and
the volume of the section of unit cube (see \cite{deboorbook},
page 2):
\begin{equation}\label{Eq:boxsplinevolume}
B(x|M)\,\,=\,\,\frac{{\rm vol}_n(P\cap
[0,1)^n)}{\sqrt{|\det(MM^T)|}},
\end{equation}
where $P:=\{y:My=x,y\in {\R}_+^n\}$. Set $A:=(a_1,\ldots,a_n)$.
Recall that $\sum_j a_j^2=1$ and $0<a_1\leq a_2\leq \cdots \leq
a_n$. Then combining (\ref{eq:def}) and (\ref{Eq:boxsplinevolume})
we have that
$$\varphi_A(\cdot)=B(\cdot+(a_1+\cdots+a_n)/2|A).$$
 Also, noting
$B(\cdot|A)$ reaching the maximum value at
$(a_1+a_2+\cdots+a_n)/2$, (\ref{eq:ba}) is equivalent to
$$
\max_x B(x |A)\,\, \geq \,\, \frac{1}{E\left(\left|\sum_{k=1}^n
a_k \varepsilon_k\right|\right)}.
$$

\section{An asymptotic formula of univariate box splines }

In this section, we shall present an asymptotic formula of
$B(\cdot|A)$. In \cite{Bconv}, Unser et. al. proved that
$B(\cdot|A)$
 tends to the Gaussian function as $n$ increase provided $a_1=a_2=\cdots=a_n$. Here,
using the saddle point approximation in statistics,  for the
general matrix $A$, we can show the box spline $B(\cdot|A)$ also
 converges   to the Gaussian function as $n$ increases:
\begin{theorem}\label{th:con}
$$
\lim_{n\rightarrow \infty}B(x|A)\,\,=\,\,
\sqrt{{6}/{\pi}}\exp(-6(x-\sum_ja_j/2)^2),
$$
where the limit may be taken pointwise or in $L^p(\R)$,
$p\in[2,\infty)$.
\end{theorem}
\begin{proof}
The saddle point approximation of $B(\cdot|A)$ is (see Theorem 6.1
and 6.2 in \cite{saddle})
\begin{equation}\label{eq:sad}
\frac{1}{(2\pi)^{1/2}|K''(s_0)|^{1/2}}\exp(K(s_0)-s_0x).
\end{equation}
Here,
$$K(s):=\ln\prod_{i}\frac{\exp(a_i s)-1}{a_i s}$$
and $s_0$ satisfies $K'(s_0)=x$.  We can consider $s_0$ as a
function of $x$. So,
\begin{equation}\label{eq:s0x}
K'(s_0)-x\,\,=\,\,0
\end{equation}
defines an implicit relationship between $s_0$ and $x$.  Noting
$K'(0)=(a_1+\cdots+a_n)/2$, the equation (\ref{eq:s0x}) implies
that $s_0((a_1+\cdots+a_n)/2)=0$. Also, by (\ref{eq:s0x}), we have
\begin{equation}\label{eq:rr}
K''(s_0)s_0'-1\,\,=\,\,0,
\end{equation}
which implies that
$$
s_0'((a_1+\cdots+a_n)/2)\,=\, {1}/{K''(0)}\,=\,
{12}/{\sum_ja_j^2}=12.
$$
Using the similar method, we have
$$
s_0''((a_1+\cdots+a_n)/2)=0,
$$
and
$$
s_0'''((a_1+\cdots+a_n)/2)=864\left(\sum_{j=1}^n
a_j^4\right)/5=O(1/n).
$$
 Using  Taylor expansion at
$(a_1+\cdots+a_n)/2$, one has
\begin{equation}\label{eq:tay}
s_0(x)\,\,=\,\,\, {12}\cdot(x-(a_1+\cdots+a_n)/2)+O(1/n).
\end{equation}
Also, by (\ref{eq:rr}), we have
$$
K''(s_0)\,\,=\,\, {1}/{s_0'}\,\,=\,\, {1}/{12}+O(1/n).
$$
 Combining Taylor expansion of $K(\cdot)$ at $0$ and
(\ref{eq:tay}), one has
$$
K(s_0)\,\,=\,\,6
{(\sum_ja_j)}(x-(a_1+\cdots+a_n)/2)+{6}(x-(a_1+\cdots+a_n)/2)^2+O(1/n).
$$
By (\ref{eq:tay}), we have
\begin{eqnarray*}
s_0\cdot x&=&{12}(x-(a_1+\cdots+a_n)/2)x+O(1/n)\\
&=&{12}(x-(a_1+\cdots+a_n)/2)^2+6(\sum_ja_j)(x-(a_1+\cdots+a_n)/2)+O(1/n).
\end{eqnarray*}
From (\ref{eq:sad}),  the saddle point approximation of $B(x|A)$
is
$$
\sqrt{{6}/{\pi}}\exp(-6(x-\sum_ja_j/2)^2)+O(1/n).
$$
The properties of the saddle point approximation imply this
theorem.
\end{proof}

\section{proof of the main result}
To prove the main result, we firstly introduce a lemma.
\begin{lemma}\label{le:1}
Put
$$
F(s)\,\, :=\,\, \frac{2}{\pi}\int_0^\infty
(1-\left|\cos(t/\sqrt{s})\right|^s)t^{-2}dt,\,\,\, s>0.
$$
Then
$$
E\left(\left|\sum_{k=1}^n a_k\varepsilon_k
\right|\right)\,\,\geq\,\, F(a_n^{-2}).
$$
\end{lemma}
\begin{proof}
By Lemma 1.3 in \cite{best}, we have
$$
E\left(|\sum_{k=1}^na_k\varepsilon_k|\right) \,\, \geq \,\,
\sum_{k=1}^n a_k^2F(a_k^{-2}).
$$
Since $F$ is an increasing  function (Lemma 1.4 \cite{best}), we
have
$$
E\left(|\sum_{k=1}^na_k\varepsilon_k|\right) \,\, \geq \,\,
\sum_{k=1}^n a_k^2F(a_k^{-2})\,\,\geq \,\,F(a_n^{-2}).
$$
\end{proof}

\begin{proof}[Proof of Theorem \ref{th:main}]
To prove the theorem, we only need to prove that there exists a
positive integer $N_0$ so that
$$
 \max_x B(x|A) \,\,\, \geq\,\,\, \frac{1}{E\left(\left|\sum_{k=1}^n
a_k \varepsilon_k\right|\right)}
$$
when $n\geq N_0$.
 Theorem \ref{th:con}
implies that
\begin{equation}\label{eq:maxlim}
\lim_{n\rightarrow \infty} \max_x B(x|A)=\sqrt{{6}/{\pi}}.
\end{equation}
Since $\sum_ja_j^2=1$ and $a_n/a_1<c_0$, we have
$\lim_{n\rightarrow\infty}1/a_n=\infty$.
 By
\cite{best}, we have $\lim_{s\rightarrow
\infty}F(s)=\sqrt{2/\pi}$. We choose $\varepsilon_0$ so that
$$0<\varepsilon_0 <\frac{\sqrt{6/\pi}-\sqrt{\pi/2}}{2}.$$
 By Lemma
\ref{le:1}, there is a positive integer $N_1$ so that
\begin{equation}\label{eq:p1}
\frac{1}{E\left(\left|\sum_{k=1}^n a_k
\varepsilon_k\right|\right)}\,\leq\, \frac{1}{F(a_n^{-2})}\,\leq\,
\sqrt{\pi/2}+\varepsilon_0
\end{equation}
provided $n\geq N_1$. The equation (\ref{eq:maxlim}) implies that
there is a positive integer $N_2(c_0)$ so that
\begin{equation}\label{eq:p2}
\max_x B(x|A)\,\,\geq\,\, \sqrt{6/\pi}-\varepsilon_0
\end{equation}
provided $n\geq N_2(c_0)$. We  set
$N_0(c_0):=\max\{N_1,N_2(c_0)\}$. Noting that
$$\sqrt{6/\pi}-\varepsilon_0\geq \sqrt{\pi/2}+\varepsilon_0,$$ combining (\ref{eq:p1}) and (\ref{eq:p2}),
 we have
$$
 \max B(x|A) \,\,\, \geq\,\,\, \frac{1}{E\left(\left|\sum_{k=1}^n
a_k \varepsilon_k\right|\right)}
$$
when $n\geq N_0(c_0)$.
\end{proof}

\bigskip

\noindent {\bf Acknowledgments.} Zhiqiang Xu was supported by
National Natural Science Foundation of China  (10871196) .
\bibliographystyle{amsplain}

\bigskip \medskip
\noindent {\bf Authors' addresses:}

\noindent Zhiqiang Xu, LSEC, Inst.~Comp.~Math., Academy of
Mathematics and System Sciences,
 Chinese Academy of Sciences, Beijing, 100080, China,
  {\tt xuzq@lsec.cc.ac.cn}

\end{document}